# On Anomaly Identification and the Counterfeit Coin Problem

Eldin Wee Chuan Lim

Department of Chemical and Biomolecular Engineering, National University of Singapore, 4 Engineering Drive 4, Singapore 117576

**Abstract**

We address a well-known problem in combinatorics involving the identification of counterfeit coins with a systematic approach. The methodology can be applied to cases where the total number of coins is exceedingly large such that brute force or enumerative comparisons become impractical. Based on the principle behind this approach, the minimum effort necessary for identification of the counterfeit coin can be determined and expressed as a simple relation. We further suggest a possible application of this methodology to error detection in quantum information processing.

## 1. INTRODUCTION

Error detection and correction techniques play an important role in the development of both classical and quantum information processing systems. The adaptation of fundamental principles and ideas of such techniques developed earlier for classical systems towards their quantum counterparts has been a popular approach in progressing research and development in the latter field in recent years. Here, we present a basic idea of a general methodology for anomaly identification in a system of discrete objects which, in our opinion, may be applicable for the development of more sophisticated error detection algorithms in various fields. We begin by illustrating this basic idea with a well-known combinatorial problem involving identification of counterfeit coins [1] and subsequently discuss the possibility of applying the general methodology to quantum error detection. We consider 12 coins which are identical in all physical aspects except for one coin which differs from the rest in weight only. A scale balance for comparing relative weights between two groups of coins is available. The required task is to identify the odd coin and conclude whether it is heavier or lighter than all other coins using the scale balance at most three times and no other equipment or approaches. Here, we will first present a systematic solution to this problem so as to illustrate the main principle behind the proposed methodology. It will then be generalized to be applicable for anomaly identification in an arbitrary system containing discrete objects, each of which may exist in any one of multiple states.

## 2. METHODOLOGY

For the present problem, the only analytical instrument to be applied is a scale which can exist in any one of three states. These correspond to situations where the group of coins placed on either side (left or right) of the scale is heavier than the other and the case where both sides exactly balance. We will use the symbols -1, 1 and 0 to represent these three states of the scale respectively. Similarly, a 3-dimensional state vector containing three state elements may be defined, such as $\begin{pmatrix} -1 \\ 1 \\ 0 \end{pmatrix}$ for example, to represent three separate analyses using the scale whereby each element denotes the state of the scale for each analysis. Further, we define the states 1 and -1 to be inverses of each other while the state 0 is the inverse of itself and two state vectors are inverses of each other if all corresponding elements of the two vectors are inverses of each other. For three applications of the scale balance, a table of all possible combinations of state vectors of the scale may be easily constructed by systematic enumeration (Table I).

**TABLE I. State vectors for a 3-state system**

|            | State Vectors |    |    |    |    |    |    |    |    |    |    |    |    |
|------------|---|---|---|---|---|---|---|---|---|---|---|---|---|
|            | 0 | 1 | 2 | 3 | 4 | 5 | 6 | 7 | 8 | 9 | 10 | 11 | 12 |
| Analysis 1 | -1 | 0 | 0 | 0 | 0 | 1 | 1 | 1 | 1 | -1 | -1 | -1 | -1 |
| Analysis 2 | -1 | 0 | -1 | -1 | -1 | 0 | 0 | 0 | 1 | -1 | 1 | 1 | 1 |
| Analysis 3 | -1 | -1 | 0 | -1 | 1 | 0 | 1 | -1 | 0 | 1 | 0 | -1 | 1 |
|            | 0' | 1' | 2' | 3' | 4' | 5' | 6' | 7' | 8' | 9' | 10' | 11' | 12' |
| Analysis 1 | 1 | 0 | 0 | 0 | 0 | -1 | -1 | -1 | -1 | 1 | 1 | 1 | 1 |
| Analysis 2 | 1 | 0 | 1 | 1 | 1 | 0 | 0 | 0 | -1 | 1 | -1 | -1 | -1 |
| Analysis 3 | 1 | 1 | 0 | 1 | -1 | 0 | -1 | 1 | 0 | -1 | 0 | 1 | -1 |

It may be noted that the state vectors (columns) presented in Table I have been arranged such that vectors 0 to 12 are inverses of vectors 0' to 12' respectively. The solution to the problem may then be completed as follows:

1. Label the twelve coins 1 to 12.
2. Place coins whose labels correspond to state elements in row 1 of Table I (named Analysis 1) which are -1 on the left side of the scale and those which are 1 on the right.
3. Record the state of the scale as the result of the first analysis.
4. Repeat using rows 2 and 3 of Table I for the second and third analyses respectively. Following this methodology, the three analyses to be conducted for the present problem and the placement of coins on each side of the scale are as shown in Table II.
5. If the resulting state vector obtained corresponds to state vector $x$ of Table I, then coin $x$ is the odd coin and it is heavier than all other coins. Conversely, if the resulting state vector corresponds to state vector $x'$, then coin $x$ is lighter than the rest.

**TABLE II. Configuration of analyses for 12-coin problem**

|            | Coin Number | | | | | | | |
|------------|---|---|---|---|---|---|---|---|
|            | Left side of scale | | | | Right side of scale | | | |
| Analysis 1 | 9 | 10 | 11 | 12 | 5 | 6 | 7 | 8 |
| Analysis 2 | 2 | 3 | 4 | 9 | 8 | 10 | 11 | 12 |
| Analysis 3 | 1 | 3 | 7 | 11 | 4 | 6 | 9 | 12 |

It may be observed from the methodology presented that every coin appears at least once in the three analyses conducted with the scale but no coin exists in all three analyses. As such, the state vectors $\begin{pmatrix}0\\0\\0\end{pmatrix}$, $\begin{pmatrix}-1\\-1\\-1\end{pmatrix}$ and $\begin{pmatrix}1\\1\\1\end{pmatrix}$ may be omitted as they represent unphysical states. We define the present system of coins and scale balance to be a 3-state system since each element of the state vector can be in any of 3 states. Based on the principle behind the above methodology, it may be easily shown that the maximum number of coins that may be present, if only three applications of the scale are allowed, is exactly 12. The number of coins which can be in the problem is also expected to be larger if the number of analyses allowed with the scale is increased. For example, with 4 analyses, a counterfeit coin may be identified from a total of 39 coins while with 5 analyses, the total number of coins increases to 120. It is straightforward to derive the following relation which relates the maximum number of coins which can be analyzed, $n_k$, to the number of analyses allowed, k.

$$n_k = \frac{1}{2}(3^k - 3) \qquad (1)$$

Alternatively, the following recurrence relations relate the maximum number of coins which can be analyzed, $n_{k+1}$, using (k + 1) analyses to that using k analyses:

$$n_{k+1} = 3^k + n_k \qquad (2)$$

$$n_{k+1} = 3(n_k + 1) \qquad (3)$$

where k is a positive integer $\geq 1$ and $n_1 = 0$.

## 3. GENERALIZATION AND POSSIBLE APPLICATION

In this section, we generalize the methodology presented above to arbitrary systems comprising many discrete objects each of which may exist in multiple states. Here, the discussion will be presented within the context of quantum information processing so as to illustrate potential applications of the proposed methodology towards quick error detection in such systems. However, it is noted that the anomaly identification methodology proposed may, in principle, also find possible applications in the biological sciences such as for error detection in gene sequencing or base pair mismatch recognition in DNA sequencing for example. We also note that some of the techniques required for practical implementation of the methodology for quantum error correction applications may not be realizable yet with present day technology. Nevertheless, it is hoped that this communication can provide a new framework for further development by research workers in this area and also inspire developers of analytical instruments to build the as yet hypothetical machines mentioned here for the further advancement of this research field.

Quantum information processing is considered to be the future of computing technology and communications. The promises of significantly larger computational power and long distance quantum cryptography for perfectly secure information transfer possible with this revolutionary technology have attracted much attention and investments from around the world. At present, some of the challenges associated with this relatively new research field include the development of algorithms for quantum system characterization, quantum error detection and quantum error correction etc. Much attention has been devoted to these and other issues by many research workers in recent years. Additionally, the notion of quantum non-demolition measurements has also grown in importance in the last decade. In particular, a number of techniques for such measurements of photon number, whereby photons whose parity is measured are not destroyed by the detectors, have been reported in the scientific literature.

We consider a quantum system consisting of a set of discrete quantum elements (such as electrons, photons or quantum dots etc) and an analytical instrument capable of non-demolition measurement of the quantum states of these elements either individually or collectively. The latter may involve Bell state measurements with many body interactions which may not be feasible yet with current technology. The state of each element has to be preserved at the end of each measurement in order for the instrument to qualify as a useful tool for error detection without interference. Recently, it has been demonstrated that non-destructive measurements of spin dynamics of an electron in a quantum dot can already be achieved using optical techniques [2]. As such, the actual realization of the methodology proposed here may already be feasible. We suppose that all elements in our quantum system were prepared to be in identical quantum states but with one exhibiting an error. As with the 12-coin problem discussed, the objective here is to identify the erroneous element using the smallest number of measurements possible so as to minimize the amount of disturbances introduced into the system. We further assume by analogy that the result of the multi-body Bell state measurement may be any of N quantum states depending on the specific states of the elements used in the measurement and their configuration. Table III shows the analogous enumeration of the state vectors for the current problem with the quantum numbers $-\frac{N}{2}$ to $\frac{N}{2}$ representing the possible states of the present N-state system. In principle, this construction would then allow the combination of elements to be used in each test with the analytical instrument to be determined. The resulting state vector from all measurements carried out would also be compared with the entries in the same table to identify deterministically the erroneous element in the system. It may be expected that such a methodology would allow much faster error identification than individual testing of each element, especially when the number of elements within the system is large, but provided that collective non-demolition measurements of quantum states are technically feasible.

**TABLE III. State vectors for a general N-state system**

| | \multicolumn{5}{c}{**State Vectors**} | | | | |
|---|---|---|---|---|---|
| | **1** | **2** | **3** | **…** | **$N^k$** |
| **Analysis 1** | $-\frac{N}{2}$ | $-\frac{N}{2}+1$ | … | … | $\frac{N}{2}$ |
| **Analysis 2** | $-\frac{N}{2}$ | $-\frac{N}{2}$ | … | … | $\frac{N}{2}$ |
| **Analysis 3** | $-\frac{N}{2}$ | $-\frac{N}{2}$ | … | … | $\frac{N}{2}$ |
| **…** | … | … | … | … | … |
| **Analysis k** | $-\frac{N}{2}$ | $-\frac{N}{2}$ | … | … | $\frac{N}{2}$ |

## 4. CONCLUSIONS

In summary, we have presented a general methodology for systematic anomaly identification in a system of discrete objects. When the number of objects is large, it may be inefficient or practically impossible to conduct individual object testing to identify the anomaly. In other cases, such as with the 12-coin problem discussed in this article, there may be constraints on the amount of testing that can be performed. The anomaly identification methodology introduced here may then be applied with apparent advantage over other ad hoc methods. To maintain generality throughout the discussion, we have not dealt with the details and intricacies associated with practical implementation of the methodology to any specific area of application but prefer to leave these to future studies. Nevertheless, we have attempted in the second part of this article to apply the methodology, in principle, to the task of quantum error detection. Here, it has been pointed out that despite the attractiveness of using this methodology to achieve quick error detection, its actual implementation in a quantum information system may call for techniques and technologies which are as yet unavailable, such as non-destructive Bell state measurements with many body interactions for example. Nevertheless, we certainly look forward to seeing the ideas presented here tested with actual physical systems to verify their potential utility in quantum information processing in the near future.

Eldin Wee Chuan LIM
National University of Singapore
Department of Chemical & Biomolecular Engineering
4 Engineering Drive 4
Singapore 117576
Email: chelwce@nus.edu.sg